\documentclass[reqno]{amsart}
\newtheorem{theorem}{Theorem}
\newtheorem{corollary}[theorem]{Corollary}
\newtheorem{lemma}[theorem]{Lemma} 

\newtheorem*{ack}{Acknowledgement}

\begin{document}
\author{Theresia Eisenk\"olbl$^1$\footnote
{$^1$ email:teisenko@radon.mat.univie.ac.at}}

\title{Rhombus Tilings of a Hexagon with Three Fixed Border Tiles}
\begin{abstract}
{We compute the number of rhombus tilings of a 
hexagon with sides $a+2,b+2,c+2,a+2,b+2,c+2$ with three fixed tiles
touching the border.
The particular case $a=b=c$ solves a problem posed by Propp.
Our result 
can also be viewed as the enumeration of plane partitions having
$a+2$ rows and $b+2$ columns, with largest entry $\le c+2$, with
a given number of entries $c+2$ in the first row, a given number of
entries 0 in the last column and a given bottom-left entry.}
\end{abstract}

\maketitle


\newfont{\fourteenpoint}{cmr10 scaled\magstep3}
\newfont{\fourteenit}{cmti10 scaled\magstep2}
\newfont{\fourteensl}{cmsl10 scaled\magstep2}
\newfont{\fourteensmc}{cmcsc10 scaled\magstep2}
\newfont{\fourteentt}{cmtt10 scaled\magstep2}
\newfont{\fourteenbf}{cmbx10 scaled\magstep2}
\newfont{\fourteeni}{cmmi10 scaled\magstep2}
\newfont{\fourteensy}{cmsy10 scaled\magstep2}
\newfont{\fourteenex}{cmex10 scaled\magstep2}
\newfont{\fourteenmsa}{msam10 scaled\magstep2}
\newfont{\fourteeneufm}{eufm10 scaled\magstep2}
\newfont{\fourteenmsb}{msbm10 scaled\magstep2}

\catcode`\@=11
\font\tenln    = line10
\font\tenlnw   = linew10

\thinlines
\newskip\Einheit \Einheit=0.6cm
\newcount\xcoord \newcount\ycoord
\newdimen\xdim \newdimen\ydim \newdimen\PfadD@cke \newdimen\Pfadd@cke
\PfadD@cke1pt \Pfadd@cke0.5pt
\def\PfadDicke#1{\PfadD@cke#1 \divide\PfadD@cke by2 \Pfadd@cke\PfadD@cke \multiply\PfadD@cke by2}
\long\def\LOOP#1\REPEAT{\def\BODY{#1}\ITERATE}
\def\ITERATE{\BODY \let\next\ITERATE \else\let\next\relax\fi \next}
\let\REPEAT=\fi
\def\Punkt{\hbox{\raise-2pt\hbox to0pt{\hss\scriptsize$\bullet$\hss}}}
\def\DuennPunkt(#1,#2){\unskip
  \raise#2 \Einheit\hbox to0pt{\hskip#1 \Einheit
          \raise-2.5pt\hbox to0pt{\hss\normalsize$\bullet$\hss}\hss}}
\def\NormalPunkt(#1,#2){\unskip
  \raise#2 \Einheit\hbox to0pt{\hskip#1 \Einheit
          \raise-3pt\hbox to0pt{\hss\large$\bullet$\hss}\hss}}
\def\DickPunkt(#1,#2){\unskip
  \raise#2 \Einheit\hbox to0pt{\hskip#1 \Einheit
          \raise-4pt\hbox to0pt{\hss\Large$\bullet$\hss}\hss}}
\def\Kreis(#1,#2){\unskip
  \raise#2 \Einheit\hbox to0pt{\hskip#1 \Einheit
          \raise-4pt\hbox to0pt{\hss\Large$\circ$\hss}\hss}}
\def\Diagonale(#1,#2)#3{\unskip\leavevmode
  \xcoord#1\relax \ycoord#2\relax
      \raise\ycoord \Einheit\hbox to0pt{\hskip\xcoord \Einheit
         \unitlength\Einheit
         \line(1,1){#3}\hss}}
\def\AntiDiagonale(#1,#2)#3{\unskip\leavevmode
  \xcoord#1\relax \ycoord#2\relax \advance\xcoord by -0.05\relax
      \raise\ycoord \Einheit\hbox to0pt{\hskip\xcoord \Einheit
         \unitlength\Einheit
         \line(1,-1){#3}\hss}}
\def\Pfad(#1,#2),#3\endPfad{\unskip\leavevmode
  \xcoord#1 \ycoord#2 \thicklines\ZeichnePfad#3\endPfad\thinlines}
\def\ZeichnePfad#1{\ifx#1\endPfad\let\next\relax
  \else\let\next\ZeichnePfad
    \ifnum#1=1
      \raise\ycoord \Einheit\hbox to0pt{\hskip\xcoord \Einheit
         \vrule height\Pfadd@cke width1 \Einheit depth\Pfadd@cke\hss}%
      \advance\xcoord by 1
    \else\ifnum#1=2
      \raise\ycoord \Einheit\hbox to0pt{\hskip\xcoord \Einheit
        \hbox{\hskip-1pt\vrule height1 \Einheit width\PfadD@cke depth0pt}\hss}%
      \advance\ycoord by 1
    \else\ifnum#1=3
      \raise\ycoord \Einheit\hbox to0pt{\hskip\xcoord \Einheit
         \unitlength\Einheit
         \line(1,1){1}\hss}
      \advance\xcoord by 1
      \advance\ycoord by 1
    \else\ifnum#1=4
      \raise\ycoord \Einheit\hbox to0pt{\hskip\xcoord \Einheit
         \unitlength\Einheit
         \line(1,-1){1}\hss}
      \advance\xcoord by 1
      \advance\ycoord by -1
    \fi\fi\fi\fi
  \fi\next}
\def\hSSchritt{\leavevmode\raise-.4pt\hbox to0pt{\hss.\hss}\hskip.2\Einheit
  \raise-.4pt\hbox to0pt{\hss.\hss}\hskip.2\Einheit
  \raise-.4pt\hbox to0pt{\hss.\hss}\hskip.2\Einheit
  \raise-.4pt\hbox to0pt{\hss.\hss}\hskip.2\Einheit
  \raise-.4pt\hbox to0pt{\hss.\hss}\hskip.2\Einheit}
\def\vSSchritt{\vbox{\baselineskip.2\Einheit\lineskiplimit0pt
\hbox{.}\hbox{.}\hbox{.}\hbox{.}\hbox{.}}}
\def\DSSchritt{\leavevmode\raise-.4pt\hbox to0pt{%
  \hbox to0pt{\hss.\hss}\hskip.2\Einheit
  \raise.2\Einheit\hbox to0pt{\hss.\hss}\hskip.2\Einheit
  \raise.4\Einheit\hbox to0pt{\hss.\hss}\hskip.2\Einheit
  \raise.6\Einheit\hbox to0pt{\hss.\hss}\hskip.2\Einheit
  \raise.8\Einheit\hbox to0pt{\hss.\hss}\hss}}
\def\dSSchritt{\leavevmode\raise-.4pt\hbox to0pt{%
  \hbox to0pt{\hss.\hss}\hskip.2\Einheit
  \raise-.2\Einheit\hbox to0pt{\hss.\hss}\hskip.2\Einheit
  \raise-.4\Einheit\hbox to0pt{\hss.\hss}\hskip.2\Einheit
  \raise-.6\Einheit\hbox to0pt{\hss.\hss}\hskip.2\Einheit
  \raise-.8\Einheit\hbox to0pt{\hss.\hss}\hss}}
\def\SPfad(#1,#2),#3\endSPfad{\unskip\leavevmode
  \xcoord#1 \ycoord#2 \ZeichneSPfad#3\endSPfad}
\def\ZeichneSPfad#1{\ifx#1\endSPfad\let\next\relax
  \else\let\next\ZeichneSPfad
    \ifnum#1=1
      \raise\ycoord \Einheit\hbox to0pt{\hskip\xcoord \Einheit
         \hSSchritt\hss}%
      \advance\xcoord by 1
    \else\ifnum#1=2
      \raise\ycoord \Einheit\hbox to0pt{\hskip\xcoord \Einheit
        \hbox{\hskip-2pt \vSSchritt}\hss}%
      \advance\ycoord by 1
    \else\ifnum#1=3
      \raise\ycoord \Einheit\hbox to0pt{\hskip\xcoord \Einheit
         \DSSchritt\hss}
      \advance\xcoord by 1
      \advance\ycoord by 1
    \else\ifnum#1=4
      \raise\ycoord \Einheit\hbox to0pt{\hskip\xcoord \Einheit
         \dSSchritt\hss}
      \advance\xcoord by 1
      \advance\ycoord by -1
    \fi\fi\fi\fi
  \fi\next}
\def\Koordinatenachsen(#1,#2){\unskip
 \hbox to0pt{\hskip-.5pt\vrule height#2 \Einheit width.5pt depth1 \Einheit}%
 \hbox to0pt{\hskip-1 \Einheit \xcoord#1 \advance\xcoord by1
    \vrule height0.25pt width\xcoord \Einheit depth0.25pt\hss}}
\def\Koordinatenachsen(#1,#2)(#3,#4){\unskip
 \hbox to0pt{\hskip-.5pt \ycoord-#4 \advance\ycoord by1
    \vrule height#2 \Einheit width.5pt depth\ycoord \Einheit}%
 \hbox to0pt{\hskip-1 \Einheit \hskip#3\Einheit 
    \xcoord#1 \advance\xcoord by1 \advance\xcoord by-#3 
    \vrule height0.25pt width\xcoord \Einheit depth0.25pt\hss}}
\def\Gitter(#1,#2){\unskip \xcoord0 \ycoord0 \leavevmode
  \LOOP\ifnum\ycoord<#2
    \loop\ifnum\xcoord<#1
      \raise\ycoord \Einheit\hbox to0pt{\hskip\xcoord \Einheit\Punkt\hss}%
      \advance\xcoord by1
    \repeat
    \xcoord0
    \advance\ycoord by1
  \REPEAT}
\def\Gitter(#1,#2)(#3,#4){\unskip \xcoord#3 \ycoord#4 \leavevmode
  \LOOP\ifnum\ycoord<#2
    \loop\ifnum\xcoord<#1
      \raise\ycoord \Einheit\hbox to0pt{\hskip\xcoord \Einheit\Punkt\hss}%
      \advance\xcoord by1
    \repeat
    \xcoord#3
    \advance\ycoord by1
  \REPEAT}
\def\Label#1#2(#3,#4){\unskip \xdim#3 \Einheit \ydim#4 \Einheit
  \def\lo{\advance\xdim by-.5 \Einheit \advance\ydim by.5 \Einheit}%
  \def\llo{\advance\xdim by-.25cm \advance\ydim by.5 \Einheit}%
  \def\loo{\advance\xdim by-.5 \Einheit \advance\ydim by.25cm}%
  \def\o{\advance\ydim by.25cm}%
  \def\ro{\advance\xdim by.5 \Einheit \advance\ydim by.5 \Einheit}%
  \def\rro{\advance\xdim by.25cm \advance\ydim by.5 \Einheit}%
  \def\roo{\advance\xdim by.5 \Einheit \advance\ydim by.25cm}%
  \def\l{\advance\xdim by-.30cm}%
  \def\r{\advance\xdim by.30cm}%
  \def\lu{\advance\xdim by-.5 \Einheit \advance\ydim by-.6 \Einheit}%
  \def\llu{\advance\xdim by-.25cm \advance\ydim by-.6 \Einheit}%
  \def\luu{\advance\xdim by-.5 \Einheit \advance\ydim by-.30cm}%
  \def\u{\advance\ydim by-.30cm}%
  \def\ru{\advance\xdim by.5 \Einheit \advance\ydim by-.6 \Einheit}%
  \def\rru{\advance\xdim by.25cm \advance\ydim by-.6 \Einheit}%
  \def\ruu{\advance\xdim by.5 \Einheit \advance\ydim by-.30cm}%
  #1\raise\ydim\hbox to0pt{\hskip\xdim
     \vbox to0pt{\vss\hbox to0pt{\hss$#2$\hss}\vss}\hss}%
}
\catcode`\@=12






\def\setRevDate $#1 #2 #3${#2}
\def\TeXdrawId{\setRevDate $Date: 1995/12/19 16:40:42 $ TeXdraw V2R0}
\chardef\catamp=\the\catcode`\@
\catcode`\@=11
\long
\def\centertexdraw #1{\hbox to \hsize{\hss
\btexdraw #1\etexdraw
\hss}}
\def\btexdraw {\x@pix=0             \y@pix=0
\x@segoffpix=\x@pix  \y@segoffpix=\y@pix
\t@exdrawdef
\setbox\t@xdbox=\vbox\bgroup\offinterlineskip
\global\d@bs=0
\global\t@extonlytrue
\global\p@osinitfalse
\s@avemove \x@pix \y@pix
\m@pendingfalse
\global\p@osinitfalse
\p@athfalse
\the\everytexdraw}
\def\etexdraw {\ift@extonly \else
\t@drclose
\fi
\egroup
\ifdim \wd\t@xdbox>0pt
\t@xderror {TeXdraw box non-zero size,
possible extraneous text}%
\fi
\vbox {\offinterlineskip
\pixtobp \xminpix \l@lxbp  \pixtobp \yminpix \l@lybp
\pixtobp \xmaxpix \u@rxbp  \pixtobp \ymaxpix \u@rybp
\hbox{\t@xdinclude 
[{\l@lxbp},{\l@lybp}][{\u@rxbp},{\u@rybp}]{\p@sfile}}%
\pixtodim \xminpix \t@xpos  \pixtodim \yminpix \t@ypos
\kern \t@ypos
\hbox {\kern -\t@xpos
\box\t@xdbox
\kern \t@xpos}%
\kern -\t@ypos\relax}}
\def\drawdim #1 {\def\d@dim{#1\relax}}
\def\setunitscale #1 {\edef\u@nitsc{#1}%
\realmult \u@nitsc \s@egsc \d@sc}
\def\relunitscale #1 {\realmult {#1}\u@nitsc \u@nitsc
\realmult \u@nitsc \s@egsc \d@sc}
\def\setsegscale #1 {\edef\s@egsc {#1}%
\realmult \u@nitsc \s@egsc \d@sc}
\def\relsegscale #1 {\realmult {#1}\s@egsc \s@egsc
\realmult \u@nitsc \s@egsc \d@sc}
\def\bsegment {\ifp@ath
\f@lushbs
\f@lushmove
\fi
\begingroup
\x@segoffpix=\x@pix
\y@segoffpix=\y@pix
\setsegscale 1
\global\advance \d@bs by 1\relax}
\def\esegment {\endgroup
\ifnum \d@bs=0
\writetx {es}%
\else
\global\advance \d@bs by -1
\fi}
\def\savecurrpos (#1 #2){\getsympos (#1 #2)\a@rgx\a@rgy
\s@etcsn \a@rgx {\the\x@pix}%
\s@etcsn \a@rgy {\the\y@pix}}
\def\savepos (#1 #2)(#3 #4){\getpos (#1 #2)\a@rgx\a@rgy
\coordtopix \a@rgx \t@pixa
\advance \t@pixa by \x@segoffpix
\coordtopix \a@rgy \t@pixb
\advance \t@pixb by \y@segoffpix
\getsympos (#3 #4)\a@rgx\a@rgy
\s@etcsn \a@rgx {\the\t@pixa}%
\s@etcsn \a@rgy {\the\t@pixb}}
\def\linewd #1 {\coordtopix {#1}\t@pixa
\f@lushbs
\writetx {\the\t@pixa\space sl}}
\def\setgray #1 {\f@lushbs
\writetx {#1 sg}}
\def\lpatt (#1){\listtopix (#1)\p@ixlist
\f@lushbs
\writetx {[\p@ixlist] sd}}
\def\lvec (#1 #2){\getpos (#1 #2)\a@rgx\a@rgy
\s@etpospix \a@rgx \a@rgy
\writeps {\the\x@pix\space \the\y@pix\space lv}}
\def\rlvec (#1 #2){\getpos (#1 #2)\a@rgx\a@rgy
\r@elpospix \a@rgx \a@rgy
\writeps {\the\x@pix\space \the\y@pix\space lv}}
\def\move (#1 #2){\getpos (#1 #2)\a@rgx\a@rgy
\s@etpospix \a@rgx \a@rgy
\s@avemove \x@pix \y@pix}
\def\rmove (#1 #2){\getpos (#1 #2)\a@rgx\a@rgy
\r@elpospix \a@rgx \a@rgy
\s@avemove \x@pix \y@pix}
\def\lcir r:#1 {\coordtopix {#1}\t@pixa
\writeps {\the\t@pixa\space cr}%
\r@elupd \t@pixa \t@pixa
\r@elupd {-\t@pixa}{-\t@pixa}}
\def\fcir f:#1 r:#2 {\coordtopix {#2}\t@pixa
\writeps {\the\t@pixa\space #1 fc}%
\r@elupd \t@pixa \t@pixa
\r@elupd {-\t@pixa}{-\t@pixa}}
\def\lellip rx:#1 ry:#2 {\coordtopix {#1}\t@pixa
\coordtopix {#2}\t@pixb
\writeps {\the\t@pixa\space \the\t@pixb\space el}%
\r@elupd \t@pixa \t@pixb
\r@elupd {-\t@pixa}{-\t@pixb}}
\def\fellip f:#1 rx:#2 ry:#3 {\coordtopix {#2}\t@pixa
\coordtopix {#3}\t@pixb
\writeps {\the\t@pixa\space \the\t@pixb\space #1 fe}%
\r@elupd \t@pixa \t@pixb
\r@elupd {-\t@pixa}{-\t@pixb}}
\def\larc r:#1 sd:#2 ed:#3 {\coordtopix {#1}\t@pixa
\writeps {\the\t@pixa\space #2 #3 ar}}
\def\ifill f:#1 {\writeps {#1 fl}}
\def\lfill f:#1 {\writeps {#1 fp}}
\def\htext #1{\def\testit {#1}%
\ifx \testit\l@paren
\let\next=\h@move
\else
\let\next=\h@text
\fi
\next {#1}}
\def\rtext td:#1 #2{\def\testit {#2}%
\ifx \testit\l@paren
\let\next=\r@move
\else
\let\next=\r@text
\fi
\next td:#1 {#2}}
\def\vtext {\rtext td:90 }
\def\textref h:#1 v:#2 {\ifx #1R%
\edef\l@stuff {\hss}\edef\r@stuff {}%
\else
\ifx #1C%
\edef\l@stuff {\hss}\edef\r@stuff {\hss}%
\else
\edef\l@stuff {}\edef\r@stuff {\hss}%
\fi
\fi
\ifx #2T%
\edef\t@stuff {}\edef\b@stuff {\vss}%
\else
\ifx #2C%
\edef\t@stuff {\vss}\edef\b@stuff {\vss}%
\else
\edef\t@stuff {\vss}\edef\b@stuff {}%
\fi
\fi}
\def\avec (#1 #2){\getpos (#1 #2)\a@rgx\a@rgy
\s@etpospix \a@rgx \a@rgy
\writeps {\the\x@pix\space \the\y@pix\space (\a@type)
\the\a@lenpix\space \the\a@widpix\space av}}
\def\ravec (#1 #2){\getpos (#1 #2)\a@rgx\a@rgy
\r@elpospix \a@rgx \a@rgy
\writeps {\the\x@pix\space \the\y@pix\space (\a@type)
\the\a@lenpix\space \the\a@widpix\space av}}
\def\arrowheadsize l:#1 w:#2 {\coordtopix{#1}\a@lenpix
\coordtopix{#2}\a@widpix}
\def\arrowheadtype t:#1 {\edef\a@type{#1}}
\def\clvec (#1 #2)(#3 #4)(#5 #6)%
{\getpos (#1 #2)\a@rgx\a@rgy
\coordtopix \a@rgx\t@pixa
\advance \t@pixa by \x@segoffpix
\coordtopix \a@rgy\t@pixb
\advance \t@pixb by \y@segoffpix
\getpos (#3 #4)\a@rgx\a@rgy
\coordtopix \a@rgx\t@pixc
\advance \t@pixc by \x@segoffpix
\coordtopix \a@rgy\t@pixd
\advance \t@pixd by \y@segoffpix
\getpos (#5 #6)\a@rgx\a@rgy
\s@etpospix \a@rgx \a@rgy
\writeps {\the\t@pixa\space \the\t@pixb\space
\the\t@pixc\space \the\t@pixd\space
\the\x@pix\space \the\y@pix\space cv}}
\def\drawbb {\bsegment
\drawdim bp
\linewd 0.24
\setunitscale {\p@sfactor}
\writeps {\the\xminpix\space \the\yminpix\space mv}%
\writeps {\the\xminpix\space \the\ymaxpix\space lv}%
\writeps {\the\xmaxpix\space \the\ymaxpix\space lv}%
\writeps {\the\xmaxpix\space \the\yminpix\space lv}%
\writeps {\the\xminpix\space \the\yminpix\space lv}%
\esegment}
\def\getpos (#1 #2)#3#4{\g@etargxy #1 #2 {} \\#3#4%
\c@heckast #3%
\ifa@st
\g@etsympix #3\t@pixa
\advance \t@pixa by -\x@segoffpix
\pixtocoord \t@pixa #3%
\fi
\c@heckast #4%
\ifa@st
\g@etsympix #4\t@pixa
\advance \t@pixa by -\y@segoffpix
\pixtocoord \t@pixa #4%
\fi}
\def\getsympos (#1 #2)#3#4{\g@etargxy #1 #2 {} \\#3#4%
\c@heckast #3%
\ifa@st \else
\t@xderror {TeXdraw: invalid symbolic coordinate}%
\fi
\c@heckast #4%
\ifa@st \else
\t@xderror {TeXdraw: invalid symbolic coordinate}%
\fi}
\def\listtopix (#1)#2{\def #2{}%
\edef\l@ist {#1 }%
\m@oretrue
\loop
\expandafter\g@etitem \l@ist \\\a@rgx\l@ist
\a@pppix \a@rgx #2%
\ifx \l@ist\empty
\m@orefalse
\fi
\ifm@ore
\repeat}
\def\realmult #1#2#3{\dimen0=#1pt
\dimen2=#2\dimen0
\edef #3{\expandafter\c@lean\the\dimen2}}
\def\intdiv #1#2#3{\t@counta=#1
\t@countb=#2
\ifnum \t@countb<0
\t@counta=-\t@counta
\t@countb=-\t@countb
\fi
\t@countd=1
\ifnum \t@counta<0
\t@counta=-\t@counta
\t@countd=-1
\fi
\t@countc=\t@counta  \divide \t@countc by \t@countb
\t@counte=\t@countc  \multiply \t@counte by \t@countb
\advance \t@counta by -\t@counte
\t@counte=-1
\loop
\advance \t@counte by 1
\ifnum \t@counte<16
\multiply \t@countc by 2
\multiply \t@counta by 2
\ifnum \t@counta<\t@countb \else
\advance \t@countc by 1
\advance \t@counta by -\t@countb
\fi
\repeat
\divide \t@countb by 2
\ifnum \t@counta<\t@countb
\advance \t@countc by 1
\fi
\ifnum \t@countd<0
\t@countc=-\t@countc
\fi
\dimen0=\t@countc sp
\edef #3{\expandafter\c@lean\the\dimen0}}
\def\coordtopix #1#2{\dimen0=#1\d@dim
\dimen2=\d@sc\dimen0
\t@counta=\dimen2
\t@countb=\s@ppix
\divide \t@countb by 2
\ifnum \t@counta<0
\advance \t@counta by -\t@countb
\else
\advance \t@counta by \t@countb
\fi
\divide \t@counta by \s@ppix
#2=\t@counta}
\def\pixtocoord #1#2{\t@counta=#1%
\multiply \t@counta by \s@ppix
\dimen0=\d@sc\d@dim
\t@countb=\dimen0
\intdiv \t@counta \t@countb #2}
\def\pixtodim #1#2{\t@countb=#1%
\multiply \t@countb by \s@ppix
#2=\t@countb sp\relax}
\def\pixtobp #1#2{\dimen0=\p@sfactor pt
\t@counta=\dimen0
\multiply \t@counta by #1%
\ifnum \t@counta < 0
\advance \t@counta by -32768
\else
\advance \t@counta by 32768
\fi
\divide \t@counta by 65536
\edef #2{\the\t@counta}}
\newcount\t@counta    \newcount\t@countb
\newcount\t@countc    \newcount\t@countd
\newcount\t@counte
\newcount\t@pixa      \newcount\t@pixb
\newcount\t@pixc      \newcount\t@pixd
\newdimen\t@xpos      \newdimen\t@ypos
\newcount\xminpix      \newcount\xmaxpix
\newcount\yminpix      \newcount\ymaxpix
\newcount\a@lenpix     \newcount\a@widpix
\newcount\x@pix        \newcount\y@pix
\newcount\x@segoffpix  \newcount\y@segoffpix
\newcount\x@savepix    \newcount\y@savepix
\newcount\s@ppix
\newcount\d@bs
\newcount\t@xdnum
\global\t@xdnum=0
\newbox\t@xdbox
\newwrite\drawfile
\newif\ifm@pending
\newif\ifp@ath
\newif\ifa@st
\newif\ifm@ore
\newif \ift@extonly
\newif\ifp@osinit
\newtoks\everytexdraw
\def\l@paren{(}
\def\a@st{*}
\catcode`\%=12
\def\p@b {
\catcode`\%=14
\catcode`\{=12  \catcode`\}=12  \catcode`\u=1 \catcode`\v=2
\def\l@br u{v  \def\r@br u}v
\catcode `\{=1  \catcode`\}=2   \catcode`\u=11 \catcode`\v=11
{\catcode`\p=12 \catcode`\t=12
\gdef\c@lean #1pt{#1}}
\def\sppix#1/#2 {\dimen0=1#2 \s@ppix=\dimen0
\t@counta=#1%
\divide \t@counta by 2
\advance \s@ppix by \t@counta
\divide \s@ppix by #1%
\t@counta=\s@ppix
\multiply \t@counta by 65536
\advance \t@counta by 32891
\divide \t@counta by 65782
\dimen0=\t@counta sp
\edef\p@sfactor {\expandafter\c@lean\the\dimen0}}
\def\g@etargxy #1 #2 #3 #4\\#5#6{\def #5{#1}%
\ifx #5\empty
\g@etargxy #2 #3 #4 \\#5#6
\else
\def #6{#2}%
\def\next {#3}%
\ifx \next\empty \else
\t@xderror {TeXdraw: invalid coordinate}%
\fi
\fi}
\def\c@heckast #1{\expandafter
\c@heckastll #1\\}
\def\c@heckastll #1#2\\{\def\testit {#1}%
\ifx \testit\a@st
\a@sttrue
\else
\a@stfalse
\fi}
\def\g@etsympix #1#2{\expandafter
\ifx \csname #1\endcsname \relax
\t@xderror {TeXdraw: undefined symbolic coordinate}%
\fi
#2=\csname #1\endcsname}
\def\s@etcsn #1#2{\expandafter
\xdef\csname#1\endcsname {#2}}
\def\g@etitem #1 #2\\#3#4{\edef #4{#2}\edef #3{#1}}
\def\a@pppix #1#2{\edef\next {#1}%
\ifx \next\empty \else
\coordtopix {#1}\t@pixa
\ifx #2\empty
\edef #2{\the\t@pixa}%
\else
\edef #2{#2 \the\t@pixa}%
\fi
\fi}
\def\s@etpospix #1#2{\coordtopix {#1}\x@pix
\advance \x@pix by \x@segoffpix
\coordtopix {#2}\y@pix
\advance \y@pix by \y@segoffpix
\u@pdateminmax \x@pix \y@pix}
\def\r@elpospix #1#2{\coordtopix {#1}\t@pixa
\advance \x@pix by \t@pixa
\coordtopix {#2}\t@pixa
\advance \y@pix by \t@pixa
\u@pdateminmax \x@pix \y@pix}
\def\r@elupd #1#2{\t@counta=\x@pix
\advance\t@counta by #1%
\t@countb=\y@pix
\advance\t@countb by #2%
\u@pdateminmax \t@counta \t@countb}
\def\u@pdateminmax #1#2{\ifnum #1>\xmaxpix
\global\xmaxpix=#1%
\fi
\ifnum #1<\xminpix
\global\xminpix=#1%
\fi
\ifnum #2>\ymaxpix
\global\ymaxpix=#2%
\fi
\ifnum #2<\yminpix
\global\yminpix=#2%
\fi}
\def\s@avemove #1#2{\x@savepix=#1\y@savepix=#2%
\m@pendingtrue
\ifp@osinit \else
\global\p@osinittrue
\global\xminpix=\x@savepix \global\yminpix=\y@savepix
\global\xmaxpix=\x@savepix \global\ymaxpix=\y@savepix
\fi}
\def\f@lushmove {\global\p@osinittrue
\ifm@pending
\writetx {\the\x@savepix\space \the\y@savepix\space mv}%
\m@pendingfalse
\p@athfalse
\fi}
\def\f@lushbs {\loop
\ifnum \d@bs>0
\writetx {bs}%
\global\advance \d@bs by -1
\repeat}
\def\h@move #1#2 #3)#4{\move (#2 #3)%
\h@text {#4}}
\def\h@text #1{\pixtodim \x@pix \t@xpos
\pixtodim \y@pix \t@ypos
\vbox to 0pt{\normalbaselines
\t@stuff
\kern -\t@ypos
\hbox to 0pt{\l@stuff
\kern \t@xpos
\hbox {#1}%
\kern -\t@xpos
\r@stuff}%
\kern \t@ypos
\b@stuff\relax}}
\def\r@move td:#1 #2#3 #4)#5{\move (#3 #4)%
\r@text td:#1 {#5}}
\def\r@text td:#1 #2{\vbox to 0pt{\pixtodim \x@pix \t@xpos
\pixtodim \y@pix \t@ypos
\kern -\t@ypos
\hbox to 0pt{\kern \t@xpos
\rottxt {#1}{\z@sb {#2}}%
\hss}%
\vss}}
\def\z@sb #1{\vbox to 0pt{\normalbaselines
\t@stuff
\hbox to 0pt{\l@stuff \hbox {#1}\r@stuff}%
\b@stuff}}
\ifx \rotatebox\@undefined
\def\rottxt #1#2{\bgroup
#2%
\egroup}
\else
\let\rottxt=\rotatebox
\fi
\ifx \t@xderror\@undefined
\let\t@xderror=\errmessage
\fi
\def\t@exdrawdef {\sppix 300/in
\drawdim in
\edef\u@nitsc {1}%
\setsegscale 1
\arrowheadsize l:0.16 w:0.08
\arrowheadtype t:T
\textref h:L v:B }
\ifx \includegraphics\@undefined
\def\t@xdinclude [#1,#2][#3,#4]#5{%
\begingroup
\message {<#5>}%
\leavevmode
\t@counta=-#1%
\t@countb=-#2%
\setbox0=\hbox{%
\includegraphics{#5}}%
\t@ypos=#4 bp%
\advance \t@ypos by -#2 bp%
\t@xpos=#3 bp%
\advance \t@xpos by -#1 bp%
\dp0=0pt \ht0=\t@ypos  \wd0=\t@xpos
\box0%
\endgroup}
\else
\let\t@xdinclude=\includegraphics
\fi
\def\writeps #1{\f@lushbs
\f@lushmove
\p@athtrue
\writetx {#1}}
\def\writetx #1{\ift@extonly
\global\t@extonlyfalse
\t@xdpsfn \p@sfile
\t@dropen \p@sfile
\fi
\w@rps {#1}}
\def\w@rps #1{\immediate\write\drawfile {#1}}
\def\t@xdpsfn #1{%
\global\advance \t@xdnum by 1
\ifnum \t@xdnum<10
\xdef #1{\jobname.ps\the\t@xdnum}
\else
\xdef #1{\jobname.p\the\t@xdnum}
\fi
}
\def\t@dropen #1{%
\immediate\openout\drawfile=#1%
\w@rps {\p@b PS-Adobe-3.0 EPSF-3.0}%
\w@rps {\p@p BoundingBox: (atend)}%
\w@rps {\p@p Title: TeXdraw drawing: #1}%
\w@rps {\p@p Pages: 1}%
\w@rps {\p@p Creator: \TeXdrawId}%
\w@rps {\p@p CreationDate: \the\year/\the\month/\the\day}%
\w@rps {50 dict begin}%
\w@rps {/mv {stroke moveto} def}%
\w@rps {/lv {lineto} def}%
\w@rps {/st {currentpoint stroke moveto} def}%
\w@rps {/sl {st setlinewidth} def}%
\w@rps {/sd {st 0 setdash} def}%
\w@rps {/sg {st setgray} def}%
\w@rps {/bs {gsave} def /es {stroke grestore} def}%
\w@rps {/fl \l@br gsave setgray fill grestore}%
\w@rps    { currentpoint newpath moveto\r@br\space def}%
\w@rps {/fp {gsave setgray fill grestore st} def}%
\w@rps {/cv {curveto} def}%
\w@rps {/cr \l@br gsave currentpoint newpath 3 -1 roll 0 360 arc}%
\w@rps    { stroke grestore\r@br\space def}%
\w@rps {/fc \l@br gsave setgray currentpoint newpath}%
\w@rps    { 3 -1 roll 0 360 arc fill grestore\r@br\space def}%
\w@rps {/ar {gsave currentpoint newpath 5 2 roll arc stroke grestore} def}%
\w@rps {/el \l@br gsave /svm matrix currentmatrix def}%
\w@rps    { currentpoint translate scale newpath 0 0 1 0 360 arc}%
\w@rps    { svm setmatrix stroke grestore\r@br\space def}%
\w@rps {/fe \l@br gsave setgray currentpoint translate scale newpath}%
\w@rps    { 0 0 1 0 360 arc fill grestore\r@br\space def}%
\w@rps {/av \l@br /hhwid exch 2 div def /hlen exch def}%
\w@rps    { /ah exch def /tipy exch def /tipx exch def}%
\w@rps    { currentpoint /taily exch def /tailx exch def}%
\w@rps    { /dx tipx tailx sub def /dy tipy taily sub def}%
\w@rps    { /alen dx dx mul dy dy mul add sqrt def}%
\w@rps    { /blen alen hlen sub def}%
\w@rps    { gsave tailx taily translate dy dx atan rotate}%
\w@rps    { (V) ah ne {blen 0 gt {blen 0 lineto} if} {alen 0 lineto} ifelse}%
\w@rps    { stroke blen hhwid neg moveto alen 0 lineto blen hhwid lineto}%
\w@rps    { (T) ah eq {closepath} if}%
\w@rps    { (W) ah eq {gsave 1 setgray fill grestore closepath} if}%
\w@rps    { (F) ah eq {fill} {stroke} ifelse}%
\w@rps    { grestore tipx tipy moveto\r@br\space def}%
\w@rps {\p@sfactor\space \p@sfactor\space scale}%
\w@rps {1 setlinecap 1 setlinejoin}%
\w@rps {3 setlinewidth [] 0 setdash}%
\w@rps {0 0 moveto}%
}
\def\t@drclose {%
\bgroup
\w@rps {stroke end showpage}%
\w@rps {\p@p Trailer:}%
\pixtobp \xminpix \l@lxbp  \pixtobp \yminpix \l@lybp
\pixtobp \xmaxpix \u@rxbp  \pixtobp \ymaxpix \u@rybp
\w@rps {\p@p BoundingBox: \l@lxbp\space \l@lybp\space
\u@rxbp\space \u@rybp}%
\w@rps {\p@p EOF}%
\egroup
\immediate\closeout\drawfile
}
\catcode`\@=\catamp

\def\ldreieck{\bsegment
  \rlvec(0.866025403784439 .5) \rlvec(0 -1)
  \rlvec(-0.866025403784439 .5)  
  \savepos(0.866025403784439 -.5)(*ex *ey)
        \esegment
  \move(*ex *ey)
        }
\def\rdreieck{\bsegment
  \rlvec(0.866025403784439 -.5) \rlvec(-0.866025403784439 -.5)  \rlvec(0 1)
  \savepos(0 -1)(*ex *ey)
        \esegment
  \move(*ex *ey)
        }
\def\rhombus{\bsegment
  \rlvec(0.866025403784439 .5) \rlvec(0.866025403784439 -.5) 
  \rlvec(-0.866025403784439 -.5)  \rlvec(0 1)        
  \rmove(0 -1)  \rlvec(-0.866025403784439 .5) 
  \savepos(0.866025403784439 -.5)(*ex *ey)
        \esegment
  \move(*ex *ey)
        }
\def\RhombusA{\bsegment
  \rlvec(0.866025403784439 .5) \rlvec(0.866025403784439 -.5) 
  \rlvec(-0.866025403784439 -.5) \rlvec(-0.866025403784439 .5) 
  \savepos(0.866025403784439 -.5)(*ex *ey)
        \esegment
  \move(*ex *ey)
        }
\def\RhombusB{\bsegment
  \rlvec(0.866025403784439 .5) \rlvec(0 -1)
  \rlvec(-0.866025403784439 -.5) \rlvec(0 1) 
  \savepos(0 -1)(*ex *ey)
        \esegment
  \move(*ex *ey)
        }
\def\RhombusC{\bsegment
  \rlvec(0.866025403784439 -.5) \rlvec(0 -1)
  \rlvec(-0.866025403784439 .5) \rlvec(0 1) 
  \savepos(0.866025403784439 -.5)(*ex *ey)
        \esegment
  \move(*ex *ey)
        }

\def\RdA{\bsegment
  \rlvec(0.866025403784439 .5) \rlvec(0.866025403784439 -.5) \lfill f:.2
  \rlvec(-0.866025403784439 -.5) \rlvec(-0.866025403784439 .5) \lfill f:.2
  \savepos(0.866025403784439 -.5)(*ex *ey)
        \esegment
  \move(*ex *ey)
        }
\def\RdB{\bsegment
  \rlvec(0.866025403784439 .5) \rlvec(0 -1) \lfill f:.6
  \rlvec(-0.866025403784439 -.5) \rlvec(0 1) \lfill f:.6
  \savepos(0 -1)(*ex *ey)
        \esegment
  \move(*ex *ey)
        }
\def\RdC{\bsegment
  \rlvec(0.866025403784439 -.5) \rlvec(0 -1)
  \rlvec(-0.866025403784439 .5) \rlvec(0 1) 
  \savepos(0.866025403784439 -.5)(*ex *ey)
        \esegment
  \move(*ex *ey)
        }

\def\LhombusA{\bsegment
  \rlvec(0.866025403784439 .5) \rlvec(0 -1) \rmove(0 1)
\rlvec(0.866025403784439 -.5) 
  \rlvec(-0.866025403784439 -.5) \rlvec(-0.866025403784439 .5) 
  \savepos(0.866025403784439 -.5)(*ex *ey)
        \esegment
  \move(*ex *ey)
        }
\def\LhombusB{\bsegment
  \rlvec(0.866025403784439 .5) \rlvec(0 -1) \rlvec(-0.866025403784439
  .5)
\rmove(0.866025403784439 -.5) 
  \rlvec(-0.866025403784439 -.5) \rlvec(0 1) 
  \savepos(0 -1)(*ex *ey)
        \esegment
  \move(*ex *ey)
        }
\def\LhombusC{\bsegment
  \rlvec(0.866025403784439 -.5) \rlvec(-0.866025403784439 -.5)
\rmove(0.866025403784439 .5)
\rlvec(0 -1)
  \rlvec(-0.866025403784439 .5) \rlvec(0 1) 
  \savepos(0.866025403784439 -.5)(*ex *ey)
        \esegment
  \move(*ex *ey)
        }

\def\hdSchritt{\bsegment
  \lpatt(.05 .13)
  \rlvec(0.866025403784439 -.5) 
  \savepos(0.866025403784439 -.5)(*ex *ey)
        \esegment
  \move(*ex *ey)
        }
\def\vdSchritt{\bsegment
  \lpatt(.05 .13)
  \rlvec(0 -1) 
  \savepos(0 -1)(*ex *ey)
        \esegment
  \move(*ex *ey)
        }

\begin{section}{Introduction}
The interest in rhombus tilings has emerged from the enumeration of
plane partitions in a given box (which was first carried out by
MacMahon \cite{MM}). The connection comes from representing each
entry by a stack of cubes and projecting the
picture to the plane. Then the box becomes a hexagon, where opposite
sides are equal, and the cubes turn into a rhombus tiling of the hexagon 
where the rhombi consist of two equilateral triangles (cf\@.
\cite{David-Tomei}). 
The number of plane partitions contained in an $a\times b\times
c$--box was first computed by MacMahon \cite{MM} and equals 
$$\frac {\prod_{k=0}^{a-1}{k!}\prod_{k=0}^{b-1}{k!}
\prod_{k=0}^{c-1}{k!}\prod_{k=0}^{a+b+c-1}{k!}} 
{\prod_{k=0}^{a+b-1}{k!}\prod_{k=0}^{b+c-1}{k!}
\prod_{k=0}^{a+c-1}{k!}}.$$

In \cite{Propp}, Propp proposed several problems regarding
"incomplete" hexagons, i.e., hexagons, where certain triangles are
missing. In particular, Problem~3 of \cite{Propp} asks for a formula     
for the number of 
rhombus tilings of a hexagon with sides $2n,2n+3,2n,2n+3,2n,2n+3$
and angles $120^\circ $, where the
middle triangle is missing on each of the longer borders.
This turns out to be a special case of the following result (see 
Corollary~\ref{propp}): 
\begin{theorem} \label{rh}
Let $a,b,c$ be nonnegative integers.
The number of rhombus tilings of a hexagon with sides
$a+2,c+2,b+2,a+2,c+2,b+2$,
with fixed tiles in positions $r,s,t$ touching the borders
$a+2,b+2,c+2$
respectively
(see Figure~\ref{semireg}a for the exact meaning of the parameters $r,s,t$)
equals 
\begin{multline*} 
(r+1)_b (s+1)_c (t+1)_a (c+3-t)_b (a+3-r)_c (b+3-s)_a
\frac{\prod_{k=0}^{a}{k!} \prod_{k=0}^{b}{k!} 
\prod_{k=0}^{c}{k!} \prod_{k=0}^{a+b+c+2}{k!}}
{\prod_{k=0}^{b+c+2}{k!}\prod_{k=0}^{a+c+2}{k!}\prod_{k=0}^{a+b+2}{k!}}
 \\
\times\Big( (a+1)(b+1)(c+1)(a+2-r)(b+2-s)(c+2-t) 
+ (a+1)(b+1)(c+1)r s t \hskip1.5cm\\
- (a+2-r)(b+2-s)(c+2-t)r s t
+(a+1)(c+1)(b+2-s)(c+2-t)r s\\
+(b+1)(a+1)(c+2-t)(a+2-r)s t
+(c+1)(b+1)(a+2-r)(b+2-s)t r \Big),
\end{multline*}
where $(a)_n=a (a+1)\cdots (a+n-1)$.
\end{theorem}
As shown in Figure~\ref{sch}, the fixed tiles determine the
tiling along the borders they touch. Thus, we can remove three
strips of triangles and end up with a hexagon with sidelengths
$a,c+3,b,a+3,c,b+3$ and missing border triangles in positions
$r,s,t$ (see Figure~\ref{trimiss}).
The special case $a=b=c=2n$ and $r=s=t=n+1$ solves
Problem~3 of \cite{Propp}. This is stated in the following corollary.
\begin{corollary} \label{propp}
The number of rhombus tilings of a hexagon with sides
$2n,2n+3,2n,2n+3,2n,2n+3$, where the
middle triangle is missing on each of the longer borders, equals
\begin{equation} \label{spez}
\left((n+2)_{2n}\right)^6\frac {\left(\prod _{k=0} ^{2n}{k!}\right)^3 \prod
_{k=0} ^{6n+2}{k!}} {\left(\prod _{k=0}
^{4n+2}{k!}\right)^3}(n+1)^3(3n+1)(3n+2)^2.
\end{equation}
\end{corollary}
\begin{figure}
\centertexdraw{
\drawdim truecm \linewd.02
\rhombus\rhombus\rhombus\rhombus\ldreieck
\move(.866025 .5) \rhombus\rhombus\rhombus\ldreieck
\move(-.866025 -.5) \rhombus\rhombus\rhombus\rhombus\rhombus\ldreieck
\move(-1.73205 -1) \rhombus\rhombus\rhombus\rhombus\rhombus\rhombus
\move(-1.73205 -1) \rdreieck \rhombus\rhombus\rhombus\rhombus\rhombus
\move(-1.73205 -2) \rdreieck \rhombus\rhombus\rhombus\rhombus
\move(-1.73205 -3) \rdreieck \rhombus\rhombus\rhombus
\move(1.73205 0)
\rlvec(0 -1) \rlvec(0.866025403784439 .5) 
\lfill f:0
\rlvec(0 1) \rlvec(-0.866025403784439 -.5) 
\lfill f:0
\move(0 -2)
\rlvec(-.866025 .5) \rlvec(-.866025 -.5)
\lfill f:0
\rlvec(.866025 -.5) \rlvec(.866025 .5)
\lfill f:0
\move(1.73205 -3)
\rlvec(0 -1) \rlvec(.866025 -.5)
\lfill f:0
\rlvec(0 1) \rlvec(-.866025 .5)
\lfill f:0
\htext(-1.5 -5.2){$b+2$}
\rtext td:0 (4.3 -2.1) {$\left. \vbox{\vskip1.65cm}\right\}c+2$}
\rtext td:-60 (-.2 0.3) {$\left\{ \vbox{\vskip2.1cm} \right.$}
\rtext td:60 (-.6 -5.1) {$\left\{ \vbox{\vskip1.6cm} \right.$}
\rtext td:0 (-2.2 -1.6) {$t \left\{ \vbox{\vskip.6cm} \right.$}
\rtext td:-60 (1.6 -5) {$ \left. \vbox{\vskip1.1cm} \right\}$}
\rtext td:60 (3.6 -.1) {$ \left. \vbox{\vskip1.1cm} \right\} $}
\htext(1.8 -5.5){$r$} \htext(3.7 .3) {$s$} \htext(-.7 0.4){$a+2$}
\htext(-2 -7){\small a. A hexagon with fixed tiles on three borders.}
\move(9 0)
\bsegment
\drawdim truecm \linewd.02
\RdA \RdB \RdA \RdA \RdB \RdB
\move (-0.866025403784439 -.5) 
\RdB \RdA \RdA \RdB \RdB \RdA
\move(-0.866025403784439 -2.5) 
\RdB \RdA \RdB
\move(.866025 .5)
\RdA \RdB \RdA \RdA \RdB \RdB
\move(-1.73205 -1)
\RdB\RdA \RdB\RdA \RdB\RdA 
\move(2.598 .5) \RdC \RdC
\move(-1.73205 -2) \RdC
\move(-1.73205 -3) \RdC
\move(1.73205 -2) \RdC
\move(-.866025 -3.5) \RdC
\linewd.1
\move(1.73205 0)
\rlvec(0 -1) \rlvec(0.866025403784439 .5) 
\rlvec(0 1) \rlvec(-0.866025403784439 -.5) 
\move(0 -2)
\rlvec(-.866025 .5) \rlvec(-.866025 -.5)
\rlvec(.866025 -.5) \rlvec(.866025 .5)
\move(1.73205 -3)
\rlvec(0 -1) \rlvec(.866025 -.5)
\rlvec(0 1) \rlvec(-.866025 .5)
\linewd.02
\htext(-2 -7){\small b. A rhombus tiling of the hexagon corres-}
\htext(-2 -7.5){\small ponding to the plane partition from \eqref{pp}.}
\rtext td:0 (-2.9 -2.6){$c+2\left\{ \vbox{ \vskip1.6cm } \right.$} 
\rtext td:60 (-.6 -5.1) {$\left\{ \vbox{\vskip1.6cm} \right.$}
\rtext td:-60 (2.5 -4.6) {$\left. \vbox{\vskip2.1cm} \right\} $}
\htext(-1.4 -5.3){$b+2$}
\htext(2.9 -5.1) {$a+2$}
\esegment
} 
\caption{\label{semireg}}
\end{figure}

\begin{figure}
\centertexdraw{
\drawdim cm
\drawdim truecm \linewd.02
\rhombus \ldreieck \rhombus \rhombus \ldreieck
\move (-0.866025403784439 -.5) 
 \rhombus \rhombus \rhombus \rhombus \rhombus
\move (-0.866025403784439 -.5) 
\rdreieck \rhombus \rhombus \rhombus \rhombus
\move(-0.866025403784439 -2.5) 
\rhombus \rhombus \ldreieck
\move(-0.866025403784439 -2.5) 
\rdreieck \rhombus \rhombus
\move (-0.866025403784439 -3.5) 
\rdreieck \rhombus
\move(1.73205 0)
\rlvec(0 -1) \rlvec(0.866025403784439 .5) 
\lfill f:0
\rlvec(0 1) \rlvec(-0.866025403784439 -.5) 
\lfill f:0
\move(0 -2)
\rlvec(-.866025 .5) \rlvec(-.866025 -.5)
\lfill f:0
\rlvec(.866025 -.5) \rlvec(.866025 .5)
\lfill f:0
\move(1.73205 -3)
\rlvec(0 -1) \rlvec(.866025 -.5)
\lfill f:0
\rlvec(0 1) \rlvec(-.866025 .5)
\lfill f:0
\lpatt(.05 .13)
\move(.866025 .5) \RhombusA
\move(-1.73205 -1) \RhombusB
\move(2.598 .5) \RhombusC \RhombusC
\move(-1.73205 -2) \RhombusC
\move(-1.73205 -3) \RhombusC
\move(0 -5) \RhombusA
\move(.866025 -4.5) \RhombusA
\move(2.598 -3.5) \RhombusB
\move(3.464 -3) \RhombusB
\rtext td:0 (-2.2 -1.6) {$t \left\{ \vbox{\vskip.6cm} \right.$}
\rtext td:-60 (1.6 -5) {$ \left. \vbox{\vskip1.1cm} \right\}$}
\rtext td:60 (3.6 -.1) {$ \left. \vbox{\vskip1.1cm} \right\} $}
\htext(1.8 -5.5){$r$} \htext(3.7 .3) {s}
}
\centerline{\small A hexagon with fixed tiles on
three borders.}
\centerline{\small The tiling of the dotted area is 
determined by the fixed tiles.}
\caption{\label{sch}}
\end{figure}

Theorem~\ref{rh} has also an interpretation in terms of plane
partitions. However, it makes only sense, if we view the plane
partitions as planar 
arrays of nonnegative integers with nonincreasing rows and columns.
For example, the plane partition of Figure~\ref{semireg}b 
is represented by the array
\begin{equation}\label{pp}
\begin{array}{ccl} 3&2&2\\3&2&2\\2&2&0\\2&1&0. \end{array}
\end{equation}
It is easy to see that a plane partition 
contained in an $(a+2)\times (b+2) \times
(c+2)$--box is represented by an array of integers
$\le c+2$ having $a+2$ rows and $b+2$ columns.
Furthermore, the fixed tiles of Theorem~\ref{rh} correspond to the
conditions that $b+2-s$ entries 
in the first row of the plane partition are equal to $c+2$, 
$r$ entries in the last column are equal to 0 and the bottom-left entry
is $c+2-t$ (cf. Figure~\ref{semireg}).
So Theorem~\ref{rh} has the following corollary.
\begin{corollary} \label{main}
The number of plane partitions contained in an $(a+2)\times (b+2) \times
(c+2)$--box having $b+2-s$ entries equal to $c+2$ in the first row, $r$ entries
equal to 0 in the last column and $c+2-t$ as the bottom-left entry
equals 
\begin{multline*} 
(r+1)_b (s+1)_c (t+1)_a (c+3-t)_b (a+3-r)_c (b+3-s)_a
\frac{\prod_{k=0}^{a}{k!} \prod_{k=0}^{b}{k!} 
\prod_{k=0}^{c}{k!} \prod_{k=0}^{a+b+c+2}{k!}}
{\prod_{k=0}^{b+c+2}{k!}\prod_{k=0}^{a+c+2}{k!}\prod_{k=0}^{a+b+2}{k!}}
 \\
\times\Big( (a+1)(b+1)(c+1)(a+2-r)(b+2-s)(c+2-t) 
+ (a+1)(b+1)(c+1)r s t \hskip1.5cm\\
- (a+2-r)(b+2-s)(c+2-t)r s t
+(a+1)(c+1)(b+2-s)(c+2-t)r s\\
+(b+1)(a+1)(c+2-t)(a+2-r)s t
+(c+1)(b+1)(a+2-r)(b+2-s)t r \Big) .
\end{multline*}
\end{corollary}

\begin{figure}
\centertexdraw{
\drawdim truecm \linewd.02
\rhombus \ldreieck \rhombus \rhombus \ldreieck
\move (-0.866025403784439 -.5) 
 \rhombus \rhombus \rhombus \rhombus \rhombus
\move (-0.866025403784439 -.5) 
\rdreieck  \rhombus \rhombus \rhombus \rhombus
\move(-0.866025403784439 -2.5) 
\rhombus \rhombus \ldreieck
\move(-0.866025403784439 -2.5) 
\rdreieck \rhombus \rhombus
\move (-0.866025403784439 -3.5) 
\rdreieck \rhombus
\rtext td:0 (4.3 -2.1){$  
    \left.\vbox{\vskip0.6cm}\right\} c $}

\rtext td:60 (-0.6 -5.1){$ \left. \vbox{\vskip0.5cm}\right\{ 
 $ }

\rtext td:60 (3.6 -1) {$ {\left. \vbox{\vskip1.1cm}  \right\}
 } $ }

\rtext td:120 (1.2 -4.8) {$ {\left. \vbox{\vskip1.1cm}  \right\{
 } $ }

\rtext td:0 (-1.4 -1.1) {$ {\left. \vbox{\vskip0.5cm} t \right\{
 } $ }

\rtext td:120 (0.1 0.1) {$ {\left. \vbox{\vskip1.1cm}  \right\}
 } $ }

\htext (-.3 0.4) {$a$}
\htext(-0.9 -5.3){$b$}
\htext(1.1 -5.1){$r$}
\htext(3.7 -.6){$s$}
}
\centerline{\small  A hexagon with sides $a,c+3,b,a+3,c,b+3$, with triangles}
\centerline{\small in positions $r,s,t$ missing, 
where $a=2$, $b=1$, $c=1$, $r=2$, $s=2$, $t=1$.}
\caption{\label{trimiss}}
\end{figure}

For the proof of Theorem~\ref{rh}, 
which we provide in Section 2, we proceed as follows. 
First, we use the fact mentioned immediately
after Theorem~\ref{rh}, that it
suffices to enumerate the rhombus tilings of a hexagon with
sides $a,c+3,b,a+3,c,b+3$ and missing border triangles in positions
$r,s,t$ as shown in Figure~\ref{trimiss}.
This can be expressed as a determinant by using the 
main theorem of nonintersecting lattice paths \cite[Cor.2]{gv}
(see also \cite[Theorem~1.2]{StemAE}). 
The determinant is then evaluated by induction, using
a determinant lemma from \cite{kratt} (see Lemma~\ref{det})
and equation~\eqref{jac}, a determinant formula 
published by Jacobi in 1841 (see \cite{ja}) but first
proved in 1819 by P.~Desnanot according to \cite{Muir}. The formula
is also closely related to C.~L.~Dodgson's condensation method.
\end{section}
\begin{section}{Proof of Theorem~\ref{rh}}
By the paragraph following Theorem~\ref{rh} it is enough to show that
the theorem holds for the
number of rhombus tilings of a hexagon with sides
$a,c+3,b,a+3,c,b+3$ and missing border triangles in positions
$r,s,t$ (see Figure~\ref{trimiss}). 

We start the proof by setting up a correspondence between these rhombus
tilings 
and certain families of nonintersecting lattice paths, where
nonintersecting means that no two paths have a common vertex. 
The reader should consult Figure~\ref{trafo} 
while reading the following passage.
Given a rhombus tiling of the region described above,
the lattice paths start on the centers of upper left diagonal edges
(the edges on the
side of length $a$) and the two extra edges parallel to it on the two
neighbouring sides. They end on the lower right edges (the edges on the 
side of length $a+3$). 
The paths are generated connecting the center of the respective edge
with the center of the edge lying opposite in the rhombus. This
process is iterated using the new edge and the second rhombus it
bounds and terminates on the lower right boundary edges. 
It is clear that paths starting at different points have no common
vertices and that an arbitrary family of nonintersecting paths from
the set of the upper left edges to the set of the lower right edges
lies completely inside the hexagon and can be converted back to a
tiling
(see Figure~\ref{trafo}b).
\begin{figure}
\centertexdraw{
 \drawdim truecm 
\linewd.08
\move(0 0)
\RhombusA \RhombusB \RhombusB \RhombusA \RhombusA 
\move (-0.866025403784439 -0.5)
\RhombusB \RhombusA \RhombusB \RhombusB \RhombusA
 \move (-0.866025403784439 -2.5)
\RhombusB \RhombusB \RhombusA
 \move (1.732050807568877 -1)
\RhombusA \RhombusB \RhombusA
\move(0.866025403784439 -2.5)
\RhombusC
\move(3.464101615137755 -1)
\RhombusC
\htext (-1.8 -6){\small a. A rhombus tiling of the hexagon of 
Figure~\ref{trimiss}.}
\htext (-.3 0.4) {$a$}
\htext(-0.9 -5.3){$b$}
\htext(3 -.2){$b+3$}
\htext(2.4 -4.6){$a+3$}
\rtext td:0 (4.4 -2.1){$ 
\left.\vbox{\vskip0.6cm}\right\} c $}
\rtext td:60 (-0.55 -5.15){$ \left. \vbox{\vskip0.6cm}\right\{ 
 $ }
\rtext td:60 (2.8 -0.5) {$ {\left. \vbox{\vskip2.1cm} \right\}
 } $ }
\rtext td:120 (2.4 -4.1) {$ {\left. \vbox{\vskip2.6cm} \right\{
 } $ }

\rtext td:0 (-2.1 -2.6) {$ {\left. \vbox{\vskip2.1cm} c+3 \right\{
 } $ }

\rtext td:120 (0.1 0.1) {$ {\left. \vbox{\vskip1.1cm} \right\}
 }$}
\move(8 0) 
\bsegment
\linewd.08
\move(0 0)
\RhombusA \RhombusB \RhombusB \RhombusA \RhombusA 
\move (-0.866025403784439 -0.5)
\RhombusB \RhombusA \RhombusB \RhombusB \RhombusA
 \move (-0.866025403784439 -2.5)
\RhombusB \RhombusB \RhombusA
 \move (1.732050807568877 -1)
\RhombusA \RhombusB \RhombusA
\move(0.866025403784439 -2.5)
\RhombusC
\move(3.464101615137755 -1)
\RhombusC

\linewd.05
\move(0.433012701892219 0.25)
\lcir r:.1
\rmove(0.866025403784439 -0.5)
\lcir r:.1
\rmove(0 -1)
\lcir r:.1
\rmove(0 -1)
\lcir r:.1
\rmove(0.866025403784439 -0.5)
\lcir r:.1
\rmove(0.866025403784439 -0.5)
\lcir r:.1

\move(2.165063509461097 -.75)
\lcir r:.1
\rmove(0.866025403784439 -0.5)
\lcir r:.1
\rmove(0 -1)
\lcir r:.1
\rmove(0.866025403784439 -0.5)
\lcir r:.1

\move(-0.433012701892219 -.25)
\lcir r:.1
\rmove(0 -1)
\lcir r:.1
\rmove(0.866025403784439 -0.5)
\lcir r:.1
\rmove(0 -1)
\lcir r:.1
\rmove(0 -1)
\lcir r:.1
\rmove(0.866025403784439 -0.5)
\lcir r:.1

\move(-0.433012701892219 -2.25)
\lcir r:.1
\rmove(0 -1)
\lcir r:.1
\rmove(0 -1)
\lcir r:.1
\rmove(0.866025403784439 -0.5)
\lcir r:.1

\linewd.05
\move(0.433012701892219 .25)
\hdSchritt \vdSchritt \vdSchritt \hdSchritt \hdSchritt
\move(-0.433012701892219 -0.25)
\vdSchritt \hdSchritt \vdSchritt \vdSchritt \hdSchritt
\move(-0.433012701892219 -2.25)
\vdSchritt \vdSchritt \hdSchritt
\move(2.165063509461097 -.75)
\hdSchritt \vdSchritt \hdSchritt
\rtext td:120 (1.2 -4.8) {$ {\left. \vbox{\vskip1.1cm}  \right\{
 } $ }

\rtext td:0 (-1.4 -1.1) {$ {\left. \vbox{\vskip0.5cm} t \right\{
 } $ }

\rtext td:60 (3.6 -1) {$ {\left. \vbox{\vskip1.1cm}  \right\}
 } $ }

\htext(1.1 -5.1){$r$}
\htext(3.7 -.6){$s$}
\esegment
\htext (7 -6){ \small b. The corresponding family of paths.}
}

\vskip1cm
\vbox{
\centertexdraw{
\drawdim truecm
\linewd.05
\move(0.433012701892219 0.25)
\lcir r:.1
\rmove(0.866025403784439 -0.5)
\lcir r:.1
\rmove(0 -1)
\lcir r:.1
\rmove(0 -1)
\lcir r:.1
\rmove(0.866025403784439 -0.5)
\lcir r:.1
\rmove(0.866025403784439 -0.5)
\lcir r:.1

\move(2.165063509461097 -.75)
\lcir r:.1
\rmove(0.866025403784439 -0.5)
\lcir r:.1
\rmove(0 -1)
\lcir r:.1
\rmove(0.866025403784439 -0.5)
\lcir r:.1

\move(-0.433012701892219 -.25)
\lcir r:.1
\rmove(0 -1)
\lcir r:.1
\rmove(0.866025403784439 -0.5)
\lcir r:.1
\rmove(0 -1)
\lcir r:.1
\rmove(0 -1)
\lcir r:.1
\rmove(0.866025403784439 -0.5)
\lcir r:.1

\move(-0.433012701892219 -2.25)
\lcir r:.1
\rmove(0 -1)
\lcir r:.1
\rmove(0 -1)
\lcir r:.1
\rmove(0.866025403784439 -0.5)
\lcir r:.1
\linewd.05
\move(0.433012701892219 .25)
\hdSchritt \vdSchritt \vdSchritt \hdSchritt \hdSchritt
\move(-0.433012701892219 -0.25)
\vdSchritt \hdSchritt \vdSchritt \vdSchritt \hdSchritt
\move(-0.433012701892219 -2.25)
\vdSchritt \vdSchritt \hdSchritt
\move(2.165063509461097 -.75)
\hdSchritt \vdSchritt \hdSchritt
\htext(13 0){}
\htext (-1.8 -7.5){\small c. The path family without the rhombi.}
\htext (6.2 -7.5){\small d. The orthogonal version of the path family.}
}
\vskip-9cm
$$
\Einheit1cm
\Gitter(6,6)(0,0)
\Koordinatenachsen(6,6)(0,0)
\Kreis(0,4.05)
\Kreis(2,1.05)
\Kreis(0,2.05)
\Kreis(1,0.05)
\Kreis(1,5.05)
\Kreis(4,3.05)
\Kreis(3,5.05)
\Kreis(5,4.05)
\Pfad(0,0),1\endPfad
\Pfad(0,0),22\endPfad
\Pfad(1,1),1\endPfad
\Pfad(1,1),22\endPfad
\Pfad(0,3),1\endPfad
\Pfad(0,3),2\endPfad
\Pfad(2,3),11\endPfad
\Pfad(2,3),22\endPfad
\Pfad(4,4),1\endPfad
\Pfad(4,4),2\endPfad
\Pfad(3,5),1\endPfad
\Pfad(1,5),1\endPfad
\Label\l{P_0}(0,2)
\Label\l{P_1}(0,4)
\Label\l{P_2}(1,5)
\Label\l{P_3}(3,5)
\Label\r{\kern6pt Q_0}(1,0)
\Label\r{\kern6pt Q_1}(2,1)
\Label\r{\kern6pt Q_2}(4,3)
\Label\r{\kern6pt Q_3}(5,4)
\hskip-2cm
$$
\vskip1.5cm
\caption{\label{trafo}}
}
\end{figure}

Then we transform the picture to orthogonal paths with positive
horizontal and negative vertical steps of unit length (see
Figure~\ref{trafo}c,d).
Let the starting points of the paths be denoted by $P_0, P_1, \dots ,
P_{a+1}$ and the endpoints by $Q_0, Q_1, \dots ,Q_{a+1}$.
Now we can easily
write down the coordinates of the starting points and the endpoints: 
\begin{alignat*}{2}
P_0 & = (0,c+2-t) ,\\
P_i & = (i-1,c+2+i), &&\text{for $i=1, \ldots ,a$,} \\
P_{a+1} &= (a+b+2-s,a+c+2) ,\\ 
Q_j &= (b+j+\chi(j \ge r),j+\chi(j \ge r) ) ,
&\quad 
&\text{for $j=0, \ldots ,a+1$} .
\end{alignat*}
Here, the symbol $\chi(j \ge r)$ equals $1$ for $j \ge r$ and 0 else.
It ensures that the missing edge on the side of length $a+3$ is 
skipped. 

Next we apply the main result for nonintersecting lattice paths 
\cite[Cor.2]{gv} (see also \cite[Theorem~1.2]{StemAE}).
This theorem says that the number of families of nonintersecting
lattice paths with path $i$ leading from $P_i$ to $Q_i$ is the
determinant of the matrix with $(i,j)$-entry the 
number of lattice paths leading from $P_i$ to $Q_j$, provided
that every two paths $P_i \to Q_j$, $P_k \to Q_l$ 
have a common vertex, if $i<j$ and $k>l$.
It is easily checked that our sets of starting and endpoints 
meet the required conditions.

Since the number of lattice paths with positive horizontal and
negative vertical steps from
$(a,b)$ to $(c,d)$ equals $\binom{c-a+b-d}{b-d}$, 
we can find the number of families of nonintersecting 
lattice paths (equivalently, the number of rhombus tilings of our 
hexagon) by 
evaluating the determinant of the matrix 
$M(a,b,c,r,s,t)= M=\left( M[i,j] \right)_{i,j=0}^{a+1}$, where

\begin{equation} \label{defM}
M[i,j]=\left \{ 
\begin{array}{cl} {\displaystyle\binom{b+c-t+2}{c+2-t-j-\chi(j \ge r)}} &
i=0, \vspace{1ex}\\ 
\vspace{1ex}
 {\displaystyle\binom{b+c+3}{b+j+\chi(j \ge r)-i+1}} & i=1, \ldots ,a
,\\
 {\displaystyle\binom{c+s}{j+\chi(j \ge r) -a-2+s}} & i=a+1.
\end{array}
 \right.
\end{equation} 

We will do this using a determinant formula due to Desnanot (see \cite{Muir}).
Given a matrix $A=(A[i,j])_{i,j=0}^{n}$, this formula states that
\begin{equation} \label{jac}
(\det A)(\det A_{0,n}^{0,n})=
(\det A_0^0)(\det A_n^n)-(\det A_0^n)(\det A_n^0) ,\end{equation} 
where $A_j^{i}$ denotes the matrix $A$ with 
row $i$ and column $j$ deleted, and $A_{0,n}^{0,n}$ denotes the
matrix $A$ with rows $0$ and $n$ and columns $0$ and $n$ deleted. 
(In general, 
given sequences of nonnegative integers $U$ and $V$, we
will use the symbol $A^U_V$ to denote the
matrix $A$ with all row indices from $U$ and all column indices 
from $V$ deleted.)

If we use \eqref{jac} for $A=M$ and $n=a+1$, we get
\begin{equation} \label{jacM}
(\det M)(\det M_{0,a+1}^{0,a+1})=
(\det M_0^0)(\det M_{a+1}^{a+1})-(\det M_0^{a+1})(\det M_{a+1}^0).\end{equation} 

In order to use \eqref{jacM} for the computation of $\det M$, 
we need to know the determinants of
$M_{0,a+1}^{0,a+1}$, $M_0^0$, $M_{a+1}^{a+1}$, $M_0^{a+1}$ and
 $M_{a+1}^0$. 
We start with the evaluation of $\det M _{0,a+1} ^{0,a+1}$. We will
employ
the following determinant lemma from \cite[Lemma~2.2]{kratt}: 

\begin{lemma} \label{det}
\begin {multline*} 
\det_{1 \le i,j \le n} \left( (X_j+A_n)(X_j+A_{n-1}) \cdots
(X_j+A_{i+1})(X_j+B_i)(X_j+B_{i-1}) \dots (X_j+B_2) \right) = \\
=\prod _{1 \le i < j \le n} ^{} {(X_i-X_j)} \prod _{2 \le i \le j
\le n} ^{}{(B_i-A_j)} .
\end{multline*}
\end{lemma}
\begin{lemma} \label{einfachematrix}
Let $1\le r \le a+1$. Then
\begin{multline*}
\det M_{0,a+1}^{0,a+1}=
\frac {((b+c+3)!)^a (a+c+2-r)! (b+r)!\prod_{1 \le i < j \le a+1}^{} {(j-i)} 
\prod_{k=0}^{a-2}{(b+c+k+4)^{a-1-k}} }
{(a+1-r)!(r-1)!\prod_{j=1}^{a+1}{\big( (b+j)!(a+c+2-j)!\big)}} .
\end{multline*}
\end{lemma}
\begin{proof}
We start by pulling out appropriate factors from columns and rows and
get 
\begin{align*} \det M_{0,a+1}^{0,a+1} 
&= \det _{1\le i,j\le a} \left( \binom{b+c+3}{b+j+\chi(j \ge r)-i+1}
\right) \\
&=(-1)^{\sum _{i=1} ^{a}{(a-i)}} \prod _{j=1}
^{a}{\frac{(b+c+3)!}{(b+j+\chi(j \ge r))!(c-j-\chi(j \ge r)+2+a)!}}\\
&\hskip2cm\times  \det _{1\le i,j\le a} \left( 
(b+j+\chi(j \ge r)-i+2)_{i-1} (j+\chi(j \ge r) -c-a-2)_{a-i} \right) .
\end{align*}
Applying Lemma~\ref{det} with $X_j=b+j+\chi(j \ge r)$, $A_k=-b-c-k-2$, 
$B_k =-k+2$  and simplifying yields the desired result. \end{proof}

The cases $r=0$ and $r\ge a+2$ can be reduced to the previous lemma
by observing that $r$ occurs only in terms of the form $\chi (j\ge
r)$. 
Since $\chi(j\ge 0) = \chi(j\ge 1)$ for $j\ge 1$ and $\chi(j\ge a+1)
= \chi(j\ge r)$ for $r\ge a+2$, $j\le a$ we have 
\begin{align} \label{rand1}
\det M_{0,a+1}^{0,a+1}(a,b,c,0,s,t)&=\det
M_{0,a+1}^{0,a+1}(a,b,c,1,s,t)\\
\label{rand2}
\det M_{0,a+1}^{0,a+1}(a,b,c,r,s,t)&=\det
M_{0,a+1}^{0,a+1}(a,b,c,a+1,s,t) \quad \text {for $r \ge a+2$.}
\end{align}

Now we express all remaining determinants in equation (\ref{jacM}) 
in terms of the determinant of 
$M _{0} ^{0}$. Whenever a matrix does not depend on some 
parameter because of a deleted row, we will use a star in place of the
parameter. It is easily checked that by appropriate relabelling of rows
and columns 
\begin{align}
\label{align:1} \det M _{a+1} ^{0}(a,b,c,r,s,*) & 
=\det M _{0} ^{0}(a,b-1,c+1,r+1,s-1,*)
,\\
\label{align:2} \det M _{a+1} ^{a+1}(a,b,c,r,*,t) & 
=\det M _{0} ^{0} (a,c,b,a+2-r,c+2-t,*)
,\\
\label{align:3} \det M _{0} ^{a+1}(a,b,c,r,*,t) & 
=\det M _{0}
^{0}(a,c-1,b+1,a+3-r,c+1-t,*).
\end{align}

The remaining task is to evaluate $\det M _{0} ^{0}$. 
We state the result for $\det M _{0} ^{0}$ in the following lemma.

\begin{lemma} \label{lastlemma}
Let $1 \le r \le a+2$. Then
\begin{multline*}  \det M _{0} ^{0}=  
\prod_{i=1} ^{a} {\left(\frac {(b+i+3)_{c+1-i}} {
(c+i+2)!} \right)}  
\frac{(s+1)_{c}}{(a+c+2)!} 
\frac { \prod_{i=1} ^{a} {i!}} {(r-1)!  (a+2-r)!}   \\
\times (c+2)_{a+1-r} (b+3)_{r-2} (c+1)_{a+2} (c+3)_a (b+3-s)_a \\
\times  \prod _{k=4} ^{a+2} {\left( (b+c+k)^{a+3-k} \right) }
\left( (b+2)(a+1)-(r-1)(b+2-s) \right). \end{multline*}
\end{lemma}
\begin{proof}
Using the argumentation preceding equations \eqref{rand1} and
\eqref{rand2} we get
\begin{align} \label{rand3}
\det M _{0} ^{0}(a,b,c,0,s,t)&=\det M _{0}^{0}(a,b,c,1,s,t)\\
\label{rand4}
\det M_{0}^{0}(a,b,c,a+3,s,t)&=\det M _{0} ^{0}(a,b,c,a+2,s,t).
\end{align}
Now we can prove the claimed expression for $\det M _{0} ^{0}$ by
induction on $a$.
It is easily checked that the statement of Lemma~\ref{lastlemma} 
holds for $a=1$.
Equation~\eqref{jac} with $A=M_0^0$ gives
\begin{equation} \label{00}
(\det M _{0} ^{0})
(\det M_{0,1,a+1}^{0,1,a+1})=
(\det M_{0,1}^{0,1})(\det M_{0,a+1}^{0,a+1})-(\det M_{0,1}^{0,a+1})(\det M_{0,a+1}^{0,1}).
\end{equation}
We will express the occurring determinants in terms of the
determinants of $M _{0} ^{0}$
and $M _{0,n} ^{0,n}$ to be able to carry out the induction.
We do this by relabelling rows and columns and get: 
\begin{align}
\det M _{0,1,a+1} ^{0,1,a+1}(a,b,c,r,*,*)&=\det M _{0,a} ^{0,a}
(a-1,b,c,r-1,*,*) \label{sch1},\\
\det M _{0,1} ^{0,1}(a,b,c,r,s,*)&=\det M _{0} ^{0}
(a-1,b,c,r-1,s,*)\label{sch2} ,\\
\det M _{0,a+1} ^{0,1}(a,b,c,r,s,*)&=\det M _{0} ^{0} (a-1,b-1,c+1,r,s-1,*),
\label{sch3}\\
\det M _{0,1} ^{0,a+1}(a,b,c,r,*,*)&=\det M _{0,a+1}
^{0,a+1}(a,b+1,c-1,r-1,*,*).
\label{sch4}
\end{align}
The matrices of the form $M _{0} ^{0}$ occurring in the above
equations \eqref{sch1}--\eqref{sch4} have parameter $(a-1)$ instead of 
$a$, so we can carry out the induction
step by 
using the values for $\det M _{0,n} ^{0,n}$ derived in
Lemma~\ref{einfachematrix} and the induction hypothesis
for $\det M _{0} ^{0}$. 
If $r \ge 2$, 
cancellation of common factors in the two sides of \eqref{00} yields
the identity 
\begin{multline*}
\big( (b+2)(a+1)-(r-1)(b+2-s)
\big) (a+b+2-s) \\
=\big( (b+2)a-(r-2)(b+2-s) \big) (a+b+2) 
- \big( (b+1)a-(r-1)
(b+2-s) \big) s ,
\end{multline*}
which is easily seen to be valid.
The case $r=1$ can be done analogously using equations 
\eqref{rand1} and \eqref{rand3}.
\end{proof}

\begin{proof}[Proof of Theorem~\ref{rh}]
Now we know all terms of equation \eqref{jacM} and can evaluate
$\det M$. 
It is indeed the expression of Theorem~\ref{rh}.
For, by plugging into equation \eqref{jacM}
the claimed formula for $\det M$ 
and the expressions derived in equations 
\eqref{align:1}--\eqref{align:3} and in Lemma~\ref{lastlemma}, 
we get an equation that can be simplified by cancelling common factors. 
If $1 \le r \le a+1$ the remaining identity is
\begin{align*}
&(b+1)(c+1)\big((b+2)(a+1)-(r-1)(b+2-s)\big)
\big((c+2)(a+1)-(a+1-r)t\big)\\
&\hskip2cm-s(c+2-t)\big((a+1)(c+1)-(a+2-r)t\big)\big((a+1)(b+1)-r(b+2-s)\big)\\
=& (a+1)(b+1)(c+1)(a+2-r)(b+2-s)(c+2-t) 
+ (a+1)(b+1)(c+1)r s t 
\\
&\hskip2cm- (a+2-r)(b+2-s)(c+2-t)r s t
+(a+1)(c+1)(b+2-s)(c+2-t)r s\\
&\hskip2cm+(b+1)(a+1)(c+2-t)(a+2-r)s t
+(c+1)(b+1)(a+2-r)(b+2-s)t r ,
\end{align*}
which is easily verified. The cases $r=0$ and $r=a+2$ can be done
analogously using equations~\eqref{rand1}, \eqref{rand2},
\eqref{rand3} and \eqref{rand4}. 
Thus the proof of Theorem~\ref{rh} is complete.
\end{proof}

\begin{ack}
The author would like to thank the referee for pointing out
Cor.~\ref{main}, and for valuable comments which helped to improve the
presentation of the results.
\end{ack} 
\end{section}

\end{document}